\documentclass[12pt,a4paper]{amsart}
\usepackage[left=2cm,right=2cm,top=2cm,bottom=2cm]{geometry}
\usepackage[english]{babel}
\usepackage{amsmath}
\usepackage{amsfonts}
\usepackage{amssymb}
\usepackage{amsthm}
\usepackage{graphicx}
\usepackage{enumitem}
\usepackage{stmaryrd}
\usepackage{hyperref}
\usepackage{bbm}
\usepackage{xcolor,mathtools}
\usepackage{float}
\usepackage[numbers,sort]{natbib}

\theoremstyle{definition}
\newtheorem{defi}{Definition}
\newtheorem{ex}[defi]{Example}
\newtheorem{rem}[defi]{Remark}
\theoremstyle{plain}
\newtheorem{teo}[defi]{Theorem}
\newtheorem{pro}[defi]{Proposition}

\renewcommand{\H}{\mathbb{H}}

\newcommand{\N}{\mathbb{N}}

\newcommand{\R}{\mathbb{R}}
\newcommand{\s}{\mathcal{S}}
\newcommand{\D}{\mathcal{D}}
\newcommand{\pp}{\partial}
\newcommand{\spa}{\operatorname{span}}

\newcommand{\se}{\subseteq}
\newcommand{\ceq}{\coloneqq}
\newcommand{\de}{\delta}
\newcommand{\la}{\lambda}

\begin{document}
\title{Submanifolds with boundary in sub-Riemannian Heisenberg Groups}



\author{Marco Di Marco}
\address[Di Marco]{ETH Z\"urich, Department of Mathematics, R\"amistrasse 101, 8092 Z\"urich, Switzerland. \newline \indent 
Dipartimento di Matematica ``T.~Levi-Civita'', Università di Padova, via Trieste 63, 35121 Padova, Italy. 
}
\email{mdimarco@ethz.ch}

\author{Davide Vittone}
\address[Vittone]{Dipartimento di Matematica ``T. Levi-Civita'', via Trieste 63, 35121 Padova, Italy.}
\email{davide.vittone@unipd.it}

\thanks{M.~D.~M. and D.~V. are supported by University of Padova and  GNAMPA of INdAM. M.~D.~M. acknowledges the support of the Swiss National Science Foundation (SNSF) Starting Grant \emph{Challenges and Breakthroughs in the Mathematics of Plasmas}, TMSGI2\textunderscore226018. D.~V. is also supported by INdAM project {\em VAC\&GMT} and by PRIN 2022PJ9EFL project {\em Geometric Measure Theory: Structure of Singular Measures, Regularity Theory and Applications in the Calculus of Variations} funded by the European Union - Next Generation EU, Mission 4, component 2 - CUP:E53D23005860006.}

\begin{abstract}
We discuss the notion of submanifolds with boundary with intrinsic $C^1$ regularity in  sub-Riemannian Heisenberg groups and we provide some examples. Eventually, we present a Stokes' Theorem for such submanifolds involving the integration of Rumin's differential forms  in Heisenberg groups.
\end{abstract}

\maketitle

\section{Introduction}
In this note we summarize some of the results of \cite{stokes} concerning, first, the notion of submanifolds with boundary with intrinsic $C^1$ regularity ($C^1_\H$-regular) in sub-Riemannian Heisenberg groups and, second,  the validity of a Stokes-type theorem for the integration of Heisenberg differential forms in the Rumin's complex. One of our goals is also that of adding to the current literature some examples of $C^1_\H$-regular submanifolds with boundary in Heisenberg groups.

Recall that, given an integer $n\geq 1$, the {\em Heisenberg group} $\H^n$ is the connected, simply connected, step 2 nilpotent Lie group associated with the Lie algebra linearly generated by elements $ X_1,\dots,X_n,Y_1,\dots,Y_n, T$ where the only non-trivial commutation relations are given by $[X_j,Y_j]=T$ for every $j=1,\dots,n$.
In exponential coordinates one can identify $\H^n=\R^{2n+1}=\R^n_x\times\R^n_y\times\R_t$ so that the left-invariant vector fields read as
\[
X_j=\partial_{x_j}-\frac{y_j}2\partial_t,\qquad Y_j=\partial_{y_j}+\frac{x_j}2\partial_t,\qquad T=\partial_t.
\]
For $\lambda>0$ the {\em dilations} $\de_\la(x,y,t)=(\la x,\la y,\la^2 t)$ define a one-parameter family of group isomorphisms. We fix a left-invariant, 1-homogeneous (with respect to dilations) and rotationally invariant (i.e., invariant with respect to rotations fixing the $t$-axis) distance $d$ in $\H^n$ (e.g. the Carnot-Carathéodory one) and for $m \geq 1$ we  define the spherical Hausdorff measure $\s^m$. The Hausdorff dimension of $\H^n$ is $Q:=2n+2$.

In Section~\ref{sec_subm} we introduce $C^1_\H$-regular submanifolds in Heisenberg groups with and without boundary (Definitions~\ref{def1} and~\ref{def2}) and we provide a few examples.  Observe (Remark~\ref{rem_folklore}) that it is not possible to produce smooth and $C^1_\H$-regular submanifolds with boundary with the topology of the disk. 
We begin with the most elementary examples, vertical planes (Example \ref{ex_piani}), and then move on to less trivial cases (Examples \ref{ex_gerono} and \ref{ex_toro}), which are compact $C^1_\H$-regular submanifolds with boundary obtained, respectively, by lifting and rotating an horizontal curve. In Section~\ref{sec_rumin} we introduce the Rumin's complex of differential forms in $\H^n$ and show that Rumin's forms can be integrated on  $C^1_\H$-regular submanifolds, eventually stating and discussing the proof of an Heisenberg version of Stokes' Theorem, Theorem~\ref{aux3}.\medskip

{\bf Acknowledgments.} The authors are grateful to Antoine Julia and Sebastiano Nicolussi Golo for their interest in the paper as well as for many valuable discussions. The authors would also like to thank Bruno Franchi for suggesting the argument in Remark~\ref{rem_folklore}.

\section{\texorpdfstring{$C^1_\H$}{C1H}-regular submanifolds}\label{sec_subm}
Let us first recall the definition of intrinsic $C^1$ (or $C^1_\H$) submanifold in $\H^n$, see e.g. \cite{franchi}.

\begin{defi}\label{def1}
    Let $1 \leq k \leq 2n$ be an integer. We say that $S\subseteq \H^n$ is a $C^1_\H$-regular submanifold of dimension $k$ 
\begin{itemize}
        \item when $1\leq k \leq n$, if $S$ is a $k$-dimensional submanifold with (classical Euclidean) $C^1$ regularity and $TS \subseteq \operatorname{span}(X_i,Y_i)_{1 \leq i \leq n}$,\vspace{.1cm}
        \item when $ n+1\leq k \leq 2n$, if, in a neighbourhood $U$ of every given point of $S$, there exists a continuous map  $f=(f_1,\dots, f_{2n+1-k}): U \to \R^{2n+1-k}$ such that $S\cap U=\lbrace p\in U: f(p)=0 \rbrace$ and $\nabla_\H f_1,\dots,\nabla_\H f_{2n+1-k}$ are linearly independent, where $\nabla_\H f_j \ceq (X_1f_j,\dots, Y_nf_j,0)\in\R^{2n+1}\equiv\H^n  $,\vspace{.1cm}
    \end{itemize} 
We also agree  that $S$ is a $C^1_\H$-regular submanifold of dimension $2n+1$ if $S$ is open, and that  $S$ is a $C^1_\H$-regular submanifold of dimension $0$ if $S$ is a collection of isolated points.
\end{defi}

When $k \geq n+1$, a $C^1_\H$ submanifold $S$ can  have fractal Euclidean dimension, see \cite{kirsc}; however, it still retains good  properties from the intrinsic viewpoint. In particular, its blow-up limit\footnote{Computed by left-translating the point to the group identity and then using Heisenberg dilations.} at any point is the vertical\footnote{As customary, we use the adjective {\em horizontal} for the distribution $\operatorname{span}(X_i,Y_i)_{1 \leq i \leq n}$ as well as for objects that are tangent to the horizontal distribution, while {\em vertical} refers to the direction $T$ or to objects to which $T$ is tangent. } $k$-plane $\spa( \nabla_\H f_1,\dots, \nabla_\H f_{2n+1-k} )^\perp\subset\H^n$, where $f$ is as in Definition~\ref{def1}. 
One can define the \emph{tangent $k$-vector} $t^\H_S$ as the unit multivector associated with $\spa( \nabla_\H f_1,\dots, \nabla_\H f_{2n+1-k} )^\perp$; this definition is well-posed up to a sign. We say that $S$ is {\em orientable} if it is possible to choose $t^\H_S$ continuously on the whole $S$.

The notion of $C^1_\H$-regular submanifold with boundary was introduced in~\cite{stokes}; from now on, we use the notation $\partial S\ceq\overline{S} \setminus S$ for the {\em boundary} of a submanifold $S$.

\begin{defi}\label{def2}
    Let $1\leq k\leq 2n+1$ be an integer. We say that $S \se \H^n$ is a $k$-dimensional $C^1_\H$ submanifold with boundary if 
    \begin{itemize}
        \item[(i)] $S$ is a $k$-dimensional $C^1_\H$ submanifold,
        \item[(ii)] $\pp S$ is a $(k-1)$-dimensional $C^1_\H$ submanifold,
        \item[(iii)] for every $p \in \pp S$ there exist a neighbourhood $U \ni p$ and a $k$-dimensional $C^1_\H$ submanifold $S'\subseteq U$ such that
        \[
        U \cap \overline{S} \subseteq S' \text{ and }(S' \setminus \overline{S}) \cap B(p,r) \neq  \emptyset \text{ for every }r>0.
         \]
    \end{itemize}
    When $S$ is oriented, an orientation is naturally induced on $\pp S$.
\end{defi}

Let us provide some examples of $C^1_\H$-regular submanifolds with boundary in the first Heisenberg group $\H^1$, where for simplicity we write $X,Y$ in place of $X_1,Y_1$.

\begin{ex}\label{ex_piani}
The simplest examples of $C^1_\H$-regular submanifolds with boundary is given by vertical half-planes. As an example consider the set
\[
P \ceq \lbrace (0,y,t) \in \H^1: y \in \R, t>0 \rbrace.
\]
It is immediate to see that $P$ is a $C^1_\H$ submanifold with boundary and that $\pp P= \lbrace (0,y,0) \in \H^1: y \in \R \rbrace$ is the integral line of the horizontal vector field $Y$.
\end{ex}

We now provide two examples of compact $C^1_\H$-regular submanifold with boundary. Let us first state the following observation, that is well-known to experts.

\begin{rem}\label{rem_folklore}
    A smooth surface in $\H^1$ with the topology of the (closed) disk and with horizontal boundary has necessarily  at least one characteristic\footnote{A point $p$ of a hypersurface $S\subset\H^n$ is called {\em characteristic} if $TS=\operatorname{span}(X_i,Y_i)_{1 \leq i \leq n}$. Observe that a smooth hypersurface without characteristic points is automatically  $C^1_\H$-regular.} point. 
    The existence of at least one characteristic point is due to topological reasons. Let us sketch a proof by contradiction; we thank B.~Franchi for suggesting us this argument. Assume that the (closed) disk $D$ has no characteristic points. For $p \in D$ we define $V(p) \in  T_pD$ to be the orthogonal projection of the vector field $T$ onto $T_pD$ with respect to the left invariant inner product making $X,Y,T$ orthonormal. We observe that the smooth vector field $V$ has no zeros on $D$; moreover, it can be easily seen that, for every $p\in\partial D$, the vector $V(p)$ is parallel to the outward pointing unit normal $\nu_{\partial D}(p)\in T_pD$  to $\partial D$. Therefore, possibly replacing $V$ with $-V$, $V$ is outward pointing at every point of $\partial D$ and, by the Poincaré-Hopf theorem (see for instance \cite{milnor}), the Euler characteristic of $D$ is 0,  a contradiction.

\end{rem}

Motivated by Remark~\ref{rem_folklore}, we conjecture that in $\H^1$ there exists no $C^1_\H$-regular 2-dimensional submanifold with boundary with the disk topology. Let us provide examples of compact orientable $C^1_\H$-regular submanifolds with boundary of different topologies.

\begin{ex}\label{ex_gerono}
Consider the curve $\gamma:[0,2\pi] \to \R^2$ (also known as \emph{lemniscate of Gerono}) given by
\[
\gamma(t)=(\cos(t),\sin(t)\cos(t)), \qquad t \in [0,2\pi].
\]
\begin{figure}[H]
    \centering
    \includegraphics[width=0.4\linewidth]{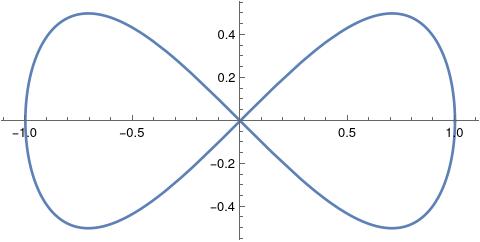}
    \caption{The lemniscate of Gerono.}
\end{figure}
Now let $\Gamma:[0,2\pi] \to \H^1$ be the \emph{horizontal lift} of $\gamma$ such that $\Gamma(0)=(1,0,0)$ (for further details on horizontal lifts see \cite[page 4]{notemonti}). The curve $\Gamma$ is a $1$-dimensional $C^1_\H$-regular submanifold (i.e., a $C^1$-regular horizontal curve) given by
\[
\Gamma(t)=\left( \cos(t),\sin(t)\cos(t),\frac{1}{24}(-9 \sin(t)-\sin(3t) )\right),\qquad t \in [0,2\pi]
\]
It is easy to check that $\Gamma$ is a closed curve without self-intersections.
\begin{figure}[H]
    \centering
    \includegraphics[width=0.3\linewidth]{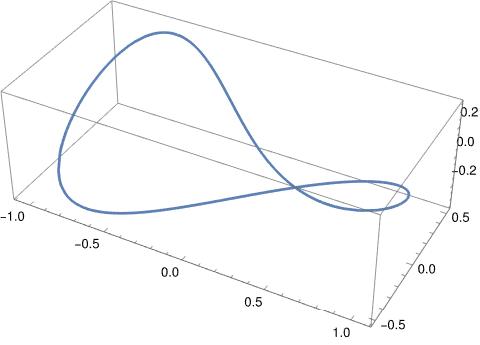}
    \includegraphics[width=0.3\linewidth]{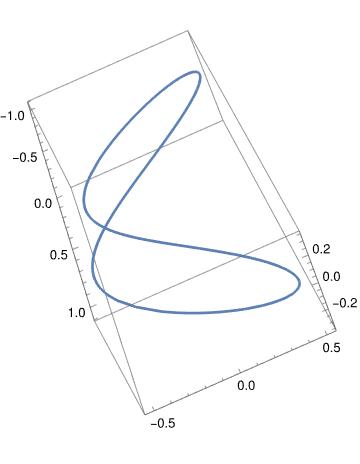}
    \includegraphics[width=0.3\linewidth]{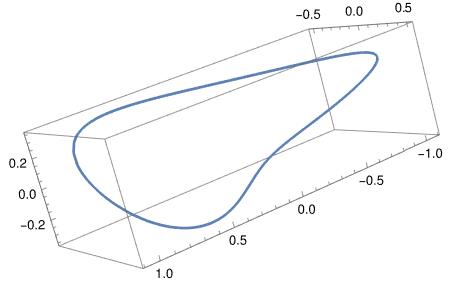}
    \caption{The horizontal lift of the lemniscate of Gerono seen from different perspectives.}
\end{figure}
Now if we consider the submanifold $\Sigma$ given from the lift (in the vertical direction) of height\footnote{We still get a set without self-intersections if we replace $1/3$ with any positive number less than $2/3$.} $1/3$, i.e.,
\[
\Sigma(t,s)=\left( \cos(t),\sin(t)\cos(t),\frac{1}{24}(-9 \sin(t)-\sin(3t) )+s\right),\qquad t \in [0,2\pi],s \in [0,1/3],
\]
we obtain an example of compact $2$-dimensional $C^1_\H$-regular submanifold with boundary. In this case, the boundary has two connected components: the horizontal curve $\Gamma$, and a vertical translation (and thus horizontal itself) of $\Gamma$.

\begin{figure}[H]
    \centering
    \includegraphics[width=0.3\linewidth]{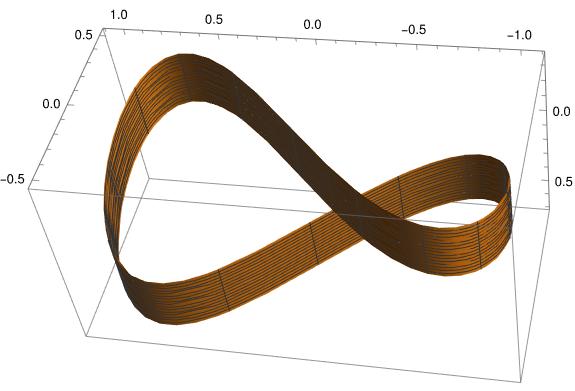}
    \includegraphics[width=0.293\linewidth]{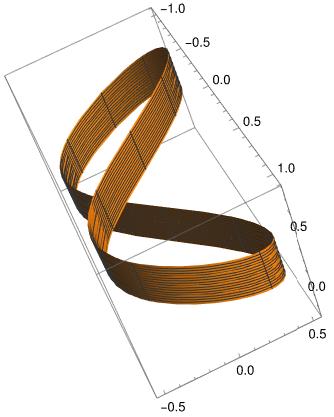}
    \includegraphics[width=0.3\linewidth]{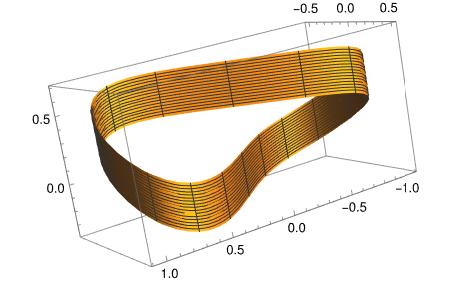}
    \caption{The compact $C^1_\H$-regular submanifold $\Sigma$ obtained by lifting the leminiscate of Gerono seen from different perspectives.}
\end{figure}
\end{ex}

\begin{ex}\label{ex_toro}
    For $R>r>0$ consider the torus $\mathcal T$ parametrized by
    \[
    \mathfrak T(u,v)= \big((R+r \cos u)\cos v, (R+r \cos u)\sin v, r \sin u\big), \qquad u,v \in [0,2\pi].
    \]
The above torus is obtained by revolving a circle of radius $r$ in the $xt$-plane around the circle of radius $R$ around the $t$-axis; $\mathcal T$ has no characteristic points. As one can see from \cite[Lemma 7.1]{bbc2022}, the characteristic foliation on an horizontal torus is filled either with periodic trajectories, or with everywhere dense trajectories,  depending only on the choice of the radii $r,R$. For example one can prove that by choosing $r=1$ and $R=\sqrt{1+n^{2/3}}$ with $n \in \N$ the characteristic foliation\footnote{The characteristic foliation is defined by the unique horizontal direction that is tangent to $\mathcal T$.} on the horizontal torus is filled with periodic trajectories (that is, closed horizontal curves) which spiral $n$ times around the torus. In Figure~\ref{fig_4} we can see one of such curves for $n=2$.

\begin{figure}[H]
    \centering
    \includegraphics[width=0.45\linewidth]{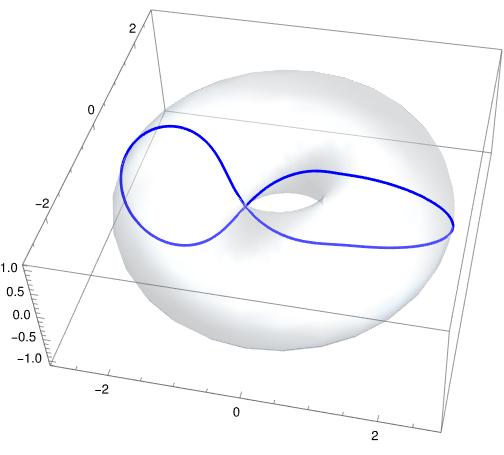}
    \includegraphics[width=0.45\linewidth]{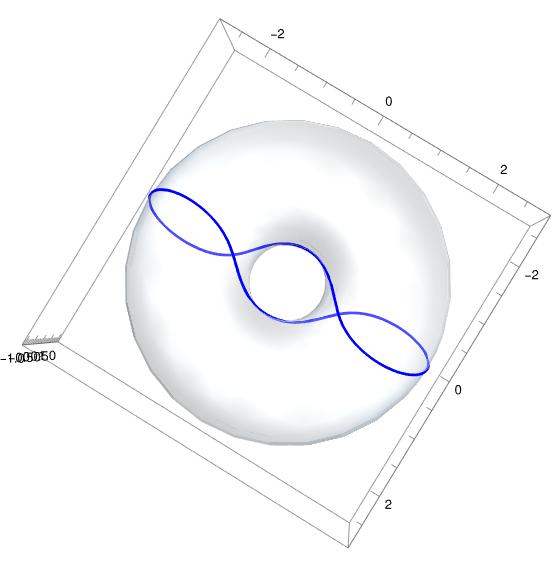}
    \caption{Two views of a closed smooth horizontal curve on $\mathcal T$ which spirals twice around the torus.}\label{fig_4}
\end{figure}

Then revolving a given periodic trajectory $\sigma$ of the characteristic foliation around the $t$-axis with an angle small enough to not have any self-intersection (in Figure \ref{fig_5} such angle is $\tfrac{\pi}{12}$), and recalling that $\spa(X,Y)$ and $\mathcal T$ are invariant under rotations around the $t$-axis, we obtain a surface $S$ on the torus (and therefore a $2$-dimensional $C^1_\H$-regular submanifold) which is bounded by two horizontal smooth curves. In conclusion $S$ is also an orientable $2$-dimensional $C^1_\H$-submanifold with boundary\footnote{As in Example \ref{ex_gerono}, the boundary of $S$ is given by two connected components: $\sigma$ and the other horizontal curve obtained by rotating $\sigma$ (with a small enough angle) around the $t$-axis.}.
\begin{figure}[H]
    \centering
    \includegraphics[width=0.45\linewidth]{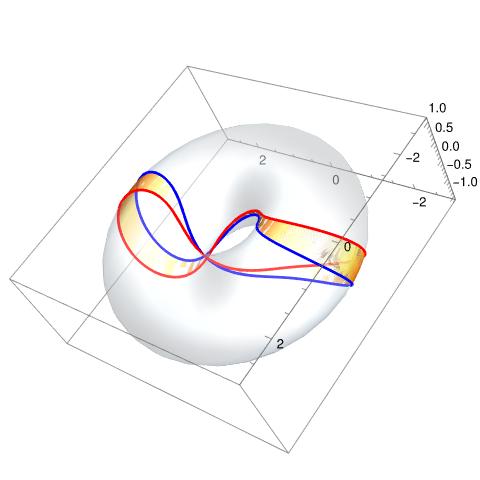}
    \includegraphics[width=0.45\linewidth]{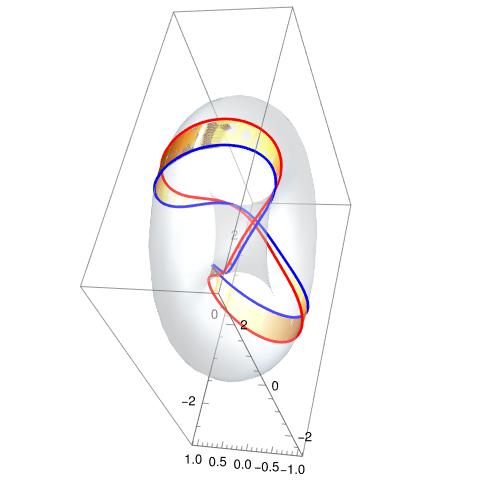}
    \caption{The compact $C^1_\H$-regular submanifold  with boundary $S$ (in yellow) obtained by revolving the curve in Figure \ref{fig_4} around the $t$-axis with an angle of $\tfrac{\pi}{12}$. The two connected components of the boundary are the original curve (in blue) and the curve rotated by $\tfrac{\pi}{12}$ (in red).}
        \label{fig_5}
\end{figure}
One can obtain more examples of $C^1_\H$-regular submanifolds with boundary starting from different curves, as Figure \ref{fig_6} shows.
\begin{figure}[H]
    \centering
    \includegraphics[width=0.45\linewidth]{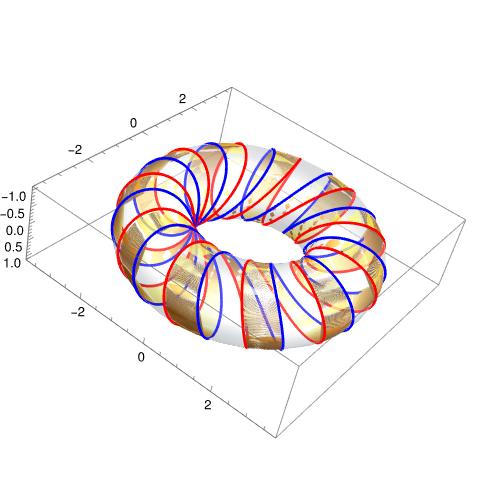}
    \includegraphics[width=0.45\linewidth]{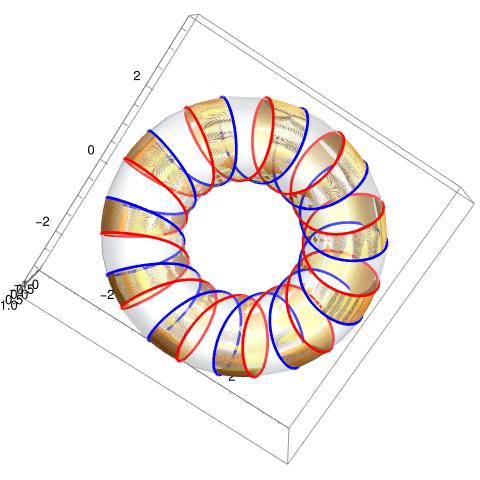}
    \caption{A compact $C^1_\H$-regular submanifold obtained as the one in Figure \ref{fig_5} starting from a curve that spirals 11 times around the torus.}
\label{fig_6}
\end{figure}

\end{ex}

\section{Rumin's complex and Stokes' Theorem in Heisenberg groups}\label{sec_rumin}

 M. Rumin introduced in \cite{RuminCR,rumin}  a complex of differential forms in $\H^n$:
    \[
    0 \xrightarrow{d} \Omega^0_\H \xrightarrow{d}  \Omega^1_\H \dots \xrightarrow{d} \Omega^n_\H \xrightarrow{D} \Omega^{n+1}_\H \xrightarrow{d} \dots \xrightarrow{d} \Omega^{2n+1}_\H \xrightarrow{d} 0.
    \]
    The are many evidences (see for instance \cite{kmx1,kmx2,kmx3,FischerTripaldi,juliapansu}) that Rumin's complex is the correct object to consider in  sub-Riemannian Heisenberg groups; let us recall that the Rumin's complex has the same cohomology as the classical De Rham's complex.
    
The precise introduction of the Rumin's complex goes beyond the scopes of this note: although the expert reader will recognize that the following presentation is an incorrect oversimplification, let us  note that
\begin{itemize}
        \item for $0\leq k\leq 2n+1$, the space $\Omega^k_\H$ of {\em Heisenberg  $k$-forms} is a linear subspace of the space of smooth differential  $k$-forms in $\H^n $,\vspace{.1cm}
        \item for $k\neq n$, $d:\Omega^k_\H\to\Omega^{k+1}_\H$ is the usual exterior differentiation operator,\vspace{.1cm}
        \item for $k=n$, the  operator $D:\Omega^n_\H\to\Omega^{n+1}_\H$ is a suitable second-order differential operator.\vspace{.1cm}
\end{itemize}
 For our scopes, it is important to notice that the structure of $D$ is
\[
D\omega=d(\omega+ \upsilon_\omega),\qquad \omega\in\Omega^n_\H
\]
where $\upsilon_\omega$ is a suitable $n$-form, depending on the first-order derivatives of $\omega$, which is ``vertical'', in the sense that it annihilates horizontal $n$-vectors.

We denote by $\D^k_\H$ the subspace of $\Omega^k_\H$ of Rumin's $k$-forms with compact support; let us introduce an integration theory for Rumin's $k$-forms  $\omega\in\D^k_\H$ on an oriented $k$-dimensional $C^1_\H$-regular submanifold $S$. This is easy when $k \leq n$, because $S$ is a submanifold with classical $C^1$ regularity and one can define ${\int_S \omega}$ by classical integration. When $k \geq n+1$, $C^1_\H$-regular submanifolds might be fractals and clearly one cannot build upon classical integration; nonetheless, when  $S$ is both $C^1$- and $C^1_\H$-regular, then (see \cite{CorniMagnani,cornimagnani2,vittone}) there exists a  constant $C=C(n,k,d)$ such that the classical integral of $\omega\in\D^{k}_\H$ on $S$ satisfies
        \[
        \int_S \omega=C \int_S \langle \omega | t^\H_S \rangle d\s^{k+1},
        \]
where $\s^h$ denotes the $h$-dimensional spherical Hausdorff measure. This suggests to define
        \[
        \int_S \omega \ceq C \int_S \langle \omega | t^\H_S \rangle d\s^{k+1}
        \]
for every oriented $C^1_\H$-regular submanifold $S$ and every $\omega\in\D^{k}_\H$.

We are finally ready to state our version of Stokes' Theorem.
\begin{teo}\label{aux3} 
Let $S \subset \H^n$ be a $k$-dimensional oriented $C^1_\H$-regular submanifold with boundary and $\omega \in \D^{k-1}_\H$. Then
\[
\int_S d_c\omega = \int_{\pp S} \omega,
\]
where $d_c=d$ if $k \neq n+1$ and $d_c=D$ if $k=n+1$.
\end{teo}

When $k \leq n$, Theorem \ref{aux3} is a direct consequence of  the classical Stokes' Theorem. Therefore, Theorem \ref{aux3} amounts to the study of the case $k \ge n+1$; here, one reasons by approximation and starts by observing that
\begin{itemize}
    \item when $k \geq n+2$  and $S$ is also  $C^1$-regular, Theorem \ref{aux3} follows directly from the classical Stokes' Theorem;\vspace{.1cm}
    \item when $k=n+1$ and $S$ is also  $C^1$-regular, Theorem \ref{aux3} follows from the classical Stokes' Theorem upon observing that 
    \[
    \int_S D\omega=\int_S d(\omega+\upsilon_\omega)=\int_{\pp S} \omega+\upsilon_\omega=\int_{\pp S}\omega,
    \]
    the last equality following from the ``verticality'' property of $\upsilon_\omega$ and the fact that the tangent space to $\pp S$ is horizontal.\vspace{.1cm}
\end{itemize}

Let us provide a sketch of the proof of Theorem \ref{aux3} by studying separately the cases $k \ge n+2$ and $k=n+1$. Details can be found in~\cite{stokes}.

\begin{proof}[Sketch of the proof of Theorem \ref{aux3}, case $k \ge n+2$]  The proof works by approximating $S$ (in $C^1_\H$) by a sequence $(S_j)_{j \in \N}$ of $C^1_\H$-and-$C^1$ submanifolds with boundary, so that
\begin{equation}\label{eq_approxxxx}
    \int_S  d_c \omega=\lim_{j \to +\infty} \int_{S_j} d_c \omega= \lim_{j \to +\infty} \int_{\pp S_j} \omega=\int_{\pp S}\omega.
\end{equation}

This can be done (locally, which is enough) exploiting the following result (see~\cite{stokes}).
    
\begin{pro}\label{prop_approx}
 Assume that $k \geq n+2$ and consider a $k$-dimensional $C^1_\H$-regular submanifold with boundary $S \subseteq \H^n$, let $h \ceq 2n+1-k$ be the codimension of $S$. Then, for every $p \in \pp S$ there exist a neighbourhood $U \ni p$ and functions $f_1,\dots, f_{h+1} \in C^1_\H(U)$ such that 
    \begin{itemize}
        \item[(i)] $\nabla_\H f_1,\dots,\nabla_\H f_{h+1}$ are linearly independent in $U$,
        \item[(ii)] $U \cap S=\lbrace q \in U: f_1(q)= \dots=f_h(q)=0,f_{h+1}(q)>0 \rbrace,$
        \item[(iii)] $U \cap \pp S=\lbrace q \in U: f_1(q)= \dots=f_h(q)=f_{h+1}(q)=0 \rbrace.$
    \end{itemize}
    Moreover, the functions $f_1,\dots,f_{h+1}$ can be chosen to be of class $C^\infty$ on $U \setminus \lbrace q \in  U : f_1(q) = \dots = f_h(q) = 0\rbrace$.
\end{pro}

The sequence $(S_j)_{j\in\N}$ needed for implementing the argument in~\eqref{eq_approxxxx} can be defined, thanks to Proposition~\ref{prop_approx}, as $S_j\ceq \lbrace q \in U: f_1(q)=1/j,\ f_2(q)= \dots=f_h(q)=0,\ f_{h+1}(q)>0 \rbrace$.
\end{proof}

Before presenting the sketch of the proof in case $k=n+1$, let us remark that Theorem \ref{aux3} was proved, although with a different formulation and only for the case $n=1$ and $k=2$, in \cite{franchitris}. The argument used in \cite{franchitris} is the following:
\begin{itemize}
\item  extend $S$ to a $C^1_\H$ surface $S'$ locally around points of $\partial S$\vspace{.1cm}
\item  $C^1_\H$-approximate $S'$ by Euclidean $C^\infty$  surfaces $S_j'$ that are non-characteristic \vspace{.1cm}
\item   produce $S_j$ by ``cutting'' $S_j'$ along an horizontal curve on $S_j'$ that is $C^1$-close to $\partial S$.\vspace{.1cm}
\end{itemize}
However, one cannot exploit this idea in $\H^n$, $n\geq 2$, because it would work only if $TS_j'\;\cap\; $ span$\{X_i,Y_i\}_{i=1,\dots,n}$ were an integrable distribution. The argument used in \cite{stokes} for every $n \geq 1$ is instead the following.

\begin{proof}[Sketch of the proof of Theorem \ref{aux3}, case $k=n+1$]  
We seek to produce a family $(S_j)_j$ of $C^1_\H$ and Euclidean $C^1$ submanifolds with boundary that approximate $S$ in $C^1_\H$, so that the strategy in~\eqref{eq_approxxxx} can be performed, and with the additional condition that $\partial S_j=\partial S$. Recall that $\partial S$ is an $n$-dimensional horizontal $C^1$ submanifold. The approximation argument goes like this:
\begin{itemize}
\item construct a $C^1_\H$ and Euclidean $C^1$  submanifold with boundary $\widetilde S$ such that $\partial \widetilde S=\partial S$. This can be done by defining $\widetilde S$ (at least locally) as the vertical half-cylinder
\[
\widetilde S:=\{(z,t):\exists (z,s)\in\pp S\text{ such that }t>s \}.
\]
Observe that such a cylinder is a $C^1_\H$ and Euclidean $C^1$ submanifold with boundary because $\pp S$ is a Euclidean $C^1$ submanifold transverse to the vertical direction $T$;\vspace{.1cm}
\item produce smooth approximations (in $C^1_\H$) $\Sigma_j$ of $S$, without any concern for their boundaries. This can be done with standard tools, e.g., by mollifying, via group convolutions, the defining function $f:\H^n\to \R^n$ of $S=\{p:f(p)=0\}$, and then considering the level sets $\Sigma_j$ of the mollified functions; such $\Sigma_j$'s are $C^1_\H$ and Euclidean smooth submanifolds that converge to $S$ in $C^1_\H$; \vspace{.1cm}
\item glue together  $\widetilde S$ (in a neighbourhood of $\partial S$) and $\Sigma_j$ (away from $\partial S$) in such a way that the resulting submanifold with boundary, $S_j$, is a $C^1_H$ and Euclidean $C^1$ submanifold with boundary that is $C^1_\H$-close to $S$.\vspace{.1cm}
\end{itemize}
This construction eventually provides the desired sequence   $(S_j)_j$. 
\end{proof}

\bibliographystyle{acm}
\bibliography{DMJNGVbib}

\end{document}